\input amstex
\input amsppt.sty
\magnification=\magstep1
\hsize=30truecc
\vsize=22.2truecm
\baselineskip=16truept
\NoBlackBoxes
\TagsOnRight \pageno=1 \nologo
\def\Z{\Bbb Z}
\def\N{\Bbb N}

\def\l{\left}
\def\r{\right}
\def\bg{\bigg}
\def\({\bg(}
\def\[{\bg\lfloor}
\def\){\bg)}
\def\]{\bg\rfloor}
\def\t{\text}
\def\f{\frac}

\def\bi{\binom}
\def\eq{\equiv}

\def\ls{\leqslant}
\def\gs{\geqslant}
\def\mo{\roman{mod}}

\def\al{\alpha}

\def\M#1#2{\thickfracwithdelims[]\thickness0{#1}{#2}_q}

\def\Proof{\noindent{\it Proof}}

\def\Remark{\medskip\noindent{\it  Remark}}

\hbox {Preprint, {\tt arXiv:0912.2671}}
\bigskip
\topmatter
\title Curious congruences for Fibonacci numbers\endtitle
\author Zhi-Wei Sun\endauthor
\leftheadtext{Zhi-Wei Sun} \rightheadtext{Curious congruences for
Fibonacci numbers}
\affil Department of Mathematics, Nanjing University\\
 Nanjing 210093, People's Republic of China
  \\  zwsun\@nju.edu.cn
  \\ {\tt http://math.nju.edu.cn/$\sim$zwsun}
\endaffil
\abstract In this paper we establish some sophisticated congruences
involving central binomial coefficients and Fibonacci numbers. For
example, we show that if $p\not=2,5$ is a prime then
$$\sum_{k=0}^{p-1}F_{2k}\bi{2k}k\eq(-1)^{\lfloor
p/5\rfloor}\l(1-\l(\f p5\r)\r)\ (\mo\ p^2)$$ and
$$\sum_{k=0}^{p-1}F_{2k+1}\bi{2k}k\eq(-1)^{\lfloor p/5\rfloor}\l(\f p5\r)\ (\mo\ p^2).$$
We also obtain similar results for some other second-order
recurrences and raise several conjectures.
\endabstract
\thanks 2010 {\it Mathematics Subject Classification}.\,Primary 11B39, 11B65;
Secondary 05A10, 11A07.
\newline\indent {\it Keywords}. Central binomial coefficients, Lucas sequences, congruences modulo prime powers.
\newline\indent Supported by the National Natural Science
Foundation (grant 10871087) and the Overseas Cooperation Fund (grant 10928101) of China.
\endthanks
\endtopmatter
\document

\heading{1. Introduction}\endheading

The well-known Fibonacci sequence $\{F_n\}_{n\gs0}$ is defined as
follows:
$$F_0=0,\ F_1=1, \t{and}\ F_{n+1}=F_n+F_{n-1}\
(n\in\Z^+=\{1,2,3,\ldots\}).$$
It plays important roles in many fields of mathematics (see, e.g., [GKP, pp.\,290-301]).

 It is known that for any odd prime $p$ we have
 $$F_p\eq\l(\f p5\r)\ (\mo\ p)\ \ \t{and}\ \ F_{p-(\f p5)}\eq0\ (\mo\ p),$$
 where $(-)$ is the Jacobi symbol. (See, e.g., [R].)

For an odd prime $p$ and an integer $m\not\eq0\ (\mo\ p)$, the sum $\sum_{k=0}^{p-1}\bi{2k}k/m^k$
and related sums modulo $p$ or $p^2$ have been investigated in [PS], [ST1], [ST2] and [S09a-S09g].

In this paper we establish some congruences involving central binomial coefficients and
Fibonacci numbers which are of a new type and
seem very curious and sophisticated.

Now we state the main results of this paper.

\proclaim{Theorem 1.1} Let $p\not=2,5$ be a prime. Then
$$\sum_{k=0}^{p-1}F_{2k}\bi{2k}k\eq(-1)^{\lfloor
p/5\rfloor}\l(1-\l(\f p5\r)\r)\ (\mo\ p^2)$$ and
$$\sum_{k=0}^{p-1}F_{2k+1}\bi{2k}k\eq(-1)^{\lfloor p/5\rfloor}\l(\f p5\r)\ (\mo\ p^2).$$
Also,
$$\sum_{k=0}^{(p-1)/2}\f{F_{2k}}{16^k}\bi{2k}k\eq(-1)^{(p-1)/2+\lfloor p/5\rfloor}\ (\mo\ p^2)$$
and
$$\sum_{k=0}^{(p-1)/2}\f{F_{2k+1}}{16^k}\bi{2k}k\eq(-1)^{(p-1)/2+\lfloor p/5\rfloor}\f{5+(\f p5)}4\ (\mo\ p^2).$$
\endproclaim
\Remark. There is no difficulty to extend Theorem 1.1 to its prime power version
(replacing $p$ in both sides of the congruences in Theorem 1.1 by $p^a$ with $a\in\Z^+$).
We can also prove the following result for any prime $p\not=2,5$ which can be viewed as a supplement to Theorem 1.1.
$$\sum_{k=0}^{p-1}F_{2k}\bi{2k}{k+1}\eq\cases0\ (\mo\ p)&\t{if}\ p\eq1\ (\mo\ 5),\\1\ (\mo\ p)&\t{if}\ p\eq-1\ (\mo\ 5),
\\-2\ (\mo\ p)&\t{if}\ p\eq-2\ (\mo\ 5),\\-3\ (\mo\ p)&\t{if}\ p\eq2\ (\mo\ 5);
\endcases$$
and
$$\sum_{k=0}^{p-1}F_{2k+1}\bi{2k}{k+1}\eq\cases0\ (\mo\ p)&\t{if}\ p\eq1\ (\mo\ 5),\\1\ (\mo\ p)&\t{if}\ p\eq-1,-2\ (\mo\ 5),
\\2\ (\mo\ p)&\t{if}\ p\eq2\ (\mo\ 5).
\endcases$$
Note that if $p$ is an odd prime and $k\in\{(p-1)/2,\ldots,p-1\}$ then $p\mid\bi{2k}k$ by Lucas' congruence (cf. [St, p.44]).

\medskip

Let $A,B\in\Z$ and $\N=\{0,1,2,\ldots\}$.  Define the Lucas sequences $u_n=u_n(A,B)\
(n\in\N)$ and $v_n=v_n(A,B)\ (n\in\N)$ as follows:
$$u_0=0,\ u_1=1,\ \t{and}\ u_{n+1}=Au_n-Bu_{n-1}\ (n=1,2,3,\ldots);$$
$$v_0=0,\ v_1=1,\ \t{and}\ v_{n+1}=Av_n-Bv_{n-1}\ (n=1,2,3,\ldots).$$
The sequence $\{u_n\}_{n\gs0}$ is a natural generalization of the Fibonacci sequence, and
$\{v_n\}_{n\gs0}$ is called the companion sequence of  $\{u_n\}_{n\gs0}$.
The characteristic equation $x^2-Ax+B=0$ of the sequences $\{u_n\}_{n\gs0}$ and $\{v_n\}_{n\gs0}$ has two roots
$$\al=\f{A+\sqrt{\Delta}}2\quad\t{and}\quad\beta=\f{A+\sqrt{\Delta}}2,$$
where $\Delta=A^2-4B$. It is well known that for any $n\in\N$ we have
$$Au_n+v_n=2u_{n+1},\ \ (\al-\beta)u_n=\al^n-\beta^n\quad\t{ and }\quad v_n=\al^n+\beta^n.$$
For convenience, we also define the sequences $\{u_n(x,y)\}_{n\gs0}$ and $\{v_n(x,y)\}_{n\gs0}$ of polynomials
as follows:
$$u_0(x,y)=0,\ u_1(x,y)=1,\ \t{and}\ u_{n+1}(x,y)=xu_n(x,y)-yu_{n-1}(x,y)\ (n\in\Z^+);$$
$$v_0(x,y)=0,\ v_1(x,y)=1,\ \t{and}\ v_{n+1}(x,y)=xv_n(x,y)-yv_{n-1}(x,y)\ (n\in\Z^+).$$

Note that $F_n=u_n(1,-1)$. Those numbers $L_n=v_n(1,-1)=2F_{n+1}-F_n$ are called Lucas numbers.
For $n\in\N$ we also have
$$\align u_n(5,5)=&\f1{\sqrt5}\(\(\f{5+\sqrt5}2\)^n-\(\f{5-\sqrt5}2\)^n\)
\\=&\sqrt5^{n-1}\(\(\f{1+\sqrt5}2\)^n-(-1)^n\(\f{1-\sqrt5}2\)^n\)
\\=&\cases 5^{n/2}F_n&\t{if}\ 2\mid n,\\5^{(n-1)/2}L_n&\t{if}\ 2\nmid n.
\endcases\endalign$$
and
$$\align v_n(5,5)=&\(\f{5+\sqrt5}2\)^n+\(\f{5-\sqrt5}2\)^n
\\=&\sqrt5^n\(\(\f{1+\sqrt5}2\)^n+(-1)^n\(\f{1-\sqrt5}2\)^n\)
\\=&\cases 5^{n/2}L_n&\t{if}\ 2\mid n,\\5^{(n+1)/2}F_n&\t{if}\ 2\nmid n.
\endcases\endalign$$

\proclaim{Theorem 1.2} Let $p\not=2,5$ be a prime. Then
$$\sum_{k=0}^{p-1}\f{u_k(5,5)}{5^k}\bi{2k}k\eq(-1)^{\lfloor p/5\rfloor}\l(\l(\f{p}5\r)-1\r)\ (\mo\ p^2)$$
and
$$\sum_{k=0}^{p-1}\f{u_{k+1}(5,5)}{5^k}\bi{2k}k\eq(-1)^{\lfloor p/5\rfloor}\l(2\l(\f{p}5\r)-1\r)\ (\mo\ p^2).$$
Also,
$$\sum_{k=0}^{(p-1)/2}\f{u_k(5,5)}{16^k}\bi{2k}k\eq\f{5(\f p5)-1}2\ (\mo\ p^2)$$
and
$$\sum_{k=0}^{(p-1)/2}\f{v_{k}(5,5)}{16^k}\bi{2k}k\eq\f{(\f p5)-1}2\ (\mo\ p^2).$$
\endproclaim

Define the sequences $\{S_n\}_{n\gs0}$ and $\{T_n\}_{n\gs0}$ as follows:
$$\align &S_0=0,\ S_1=1,\ \t{and}\ S_{n+1}=4S_n-S_{n-1}\ (n=1,2,3,\ldots);
\\&T_0=2,\ T_1=4,\ \t{and}\ T_{n+1}=4T_n-T_{n-1}\ (n=1,2,3,\ldots).
\endalign$$
Note that $S_n=u_n(4,1)$ and $T_n=v_n(4,1)$. These two sequences are also useful; see, e.g., [Sl] and [S02].

\proclaim{Theorem 1.3} Let $p>3$ be a prime. Then
$$\sum_{k=0}^{p-1}S_k\bi{2k}k\eq2\l(\l(\f p3\r)-\l(\f{-1}p\r)\r)\ (\mo\ p^2)$$
and
$$\sum_{k=0}^{p-1}S_{k+1}\bi{2k}k\eq\l(\f p3\r)\ (\mo\ p^2).$$
Also,
$$\sum_{k=0}^{(p-1)/2}\f{S_k}{16^k}\bi{2k}k\eq\f{(\f6p)-(\f 2p)}2\ (\mo\ p^2)$$
and
$$\sum_{k=0}^{(p-1)/2}\f{T_k}{16^k}\bi{2k}k\eq3\l(\f6p\r)-\l(\f 2p\r)\ (\mo\ p^2).$$
\endproclaim

The Pell sequence $\{P_n\}_{n\gs0}$ and its companion $\{Q_n\}_{n\gs0}$
are given by $P_n=u_n(2,-1)$ and $Q_n=v_n(2,1)$. For $n\in\N$ we can easily see that
$$u_n(4,2)=\cases2^{n/2}P_n&\t{if}\ 2\mid n,\\2^{(n-3)/2}Q_n&\t{if}\ 2\nmid n,\endcases
\ \t{and}\ v_n(4,2)=\cases2^{n/2}Q_n&\t{if}\ 2\mid n,\\2^{(n+3)/2}P_n&\t{if}\ 2\nmid n.\endcases
$$

\proclaim{Theorem 1.4} Let $p$ be an odd prime. Then
$$\sum_{k=0}^{p-1}\f{u_k(4,2)}{2^k}\bi{2k}k\eq\l(\f{-1}p\r)-\l(\f{-2}p\r)\ (\mo\ p^2)$$
and
$$\sum_{k=0}^{p-1}\f{u_{k+1}(4,2)}{2^k}\bi{2k}k\eq\l(\f{-1}p\r)\ (\mo\ p^2).$$
Also,
$$\sum_{k=0}^{(p-1)/2}\f{u_k(4,2)}{16^k}\bi{2k}k\eq\f{(-1)^{\lfloor (p-4)/8\rfloor}}2\l(1-\l(\f 2p\r)\r)\ (\mo\ p^2)$$
and
$$\sum_{k=0}^{(p-1)/2}\f{v_k(4,2)}{16^k}\bi{2k}k\eq2(-1)^{\lfloor p/8\rfloor}\l(\f{-1}p\r)\ (\mo\ p^2).$$
\endproclaim

We will present several lemmas in Section 2 and prove Theorems 1.1-1.4 in Section 3.
The last section contains several open conjectures.

A key point in our proofs is the use of Chebyshev polynomials. Recall that the Chebyshev polynomials
of the second kind are given by
$$U_n(\cos \theta)=\f{\sin(n+1)\theta}{\sin\theta}\ \ (n=0,1,2,\ldots).$$

\heading{2. Some lemmas}\endheading

\proclaim{Lemma 2.1} Let $p$ be any prime and let $\al$ be an algebraic integer.
Then
$$\sum_{k=0}^{p-1}\bi{2k}k\al^{p-1-k}\eq 2u_{p}(\al,\al)-u_{p}(\al-2,1)\ (\mo\ p^2).$$
\endproclaim
\Proof. In view of [ST2, Theorem 2.1],
$$\sum_{k=0}^{p-1}\bi{2k}k\al^{p-1-k}=\sum_{k=0}^{p-1}\bi{2p}ku_{p-k}(\al-2,1).$$
For $k\in\{1,\ldots,p-1\}$, we clearly have
$$\f{\bi{2p}k}{\bi pk}=\prod_{j=0}^{k-1}\f{2p-j}{p-j}\eq2\ (\mo\ p)$$
and hence
$$\bi{2p}k\eq2\bi{p}k\ (\mo\ p^2)$$
since $p\mid\bi pk$.
Therefore
$$\align&\sum_{k=0}^{p-1}\bi{2k}k\al^{p-1-k}+u_{p}(\al-2,1)
\\=&\sum_{k=1}^{p-1}\bi{2p}ku_{p-k}(\al-2,1)+2u_{p}(\al-2,1)
\\\eq&2\sum_{k=0}^{p}\bi{p}ku_{p-k}(\al-2,1)=2\sum_{j=0}^{p}\bi{p}ju_j(\al-2,1)\ (\mo\ p^2).
\endalign$$

Now it remains to show that
$$\sum_{j=0}^p\bi pju_j(\al-2,1)=u_p(\al,\al).$$
Recall that
$$u_k(x,1)=\f1{\sqrt{x^2-4}}\(\(\f{x+\sqrt{x^2-4}}2\)^k-\(\f{x-\sqrt{x^2-4}}2\)^k\)\quad\t{for all}\ k\in\N.$$
So we have
$$\align\sum_{k=0}^p\bi pk u_k(x,1)=&\f1{\sqrt{x^2-4}}\(\(\f{x+2+\sqrt{x^2-4}}2\)^p-\(\f{x+2-\sqrt{x^2-4}}2\)^p\)
\\=&u_p(x+2,x+2).\endalign$$
This concludes the proof. \qed

\medskip

\proclaim{Lemma 2.2} For any $n\in\N$ we have
$$u_{n+1}(x,1)=U_n\l(\f x2\r).$$
\endproclaim
\Proof. It is well known that
$$u_{n+1}(x,y)=\sum_{k=0}^{\lfloor n/2\rfloor}\bi{n-k}kx^{n-2k}(-y)^k$$
and
$$U_n(x)=\sum_{k=0}^{\lfloor n/2\rfloor}\bi{n-k}k(-1)^k(2x)^{n-2k}.$$
So the desired equality follows. \qed

\proclaim{Lemma 2.3} Let $p$ be a prime and let $\al\in\{(1\pm\sqrt5)/2\}$. If $\al=(1+\sqrt5)/2$, then
$$\sum_{k=0}^{p-1}\bi{2k}k\al^{2k}
\eq2\al^{p-1}\f{\sin(2p\pi/5)}{\sin(2\pi/5)}-\al^{2p-2}\f{\sin(p\pi/5)}{\sin(\pi/5)}\ (\mo\ p^2).$$
If $\al=(1-\sqrt5)/2$, then
$$\sum_{k=0}^{p-1}\bi{2k}k\al^{2k}
\eq2\al^{p-1}\f{\sin(p\pi/5)}{\sin(\pi/5)}-\al^{2p-2}\f{\sin(2p\pi/5)}{\sin(2\pi/5)}\ (\mo\ p^2).$$
\endproclaim
\Proof. Set $\beta=1-\al=-\al^{-1}$. Then
$\{\al,\beta\}=\{(1\pm \sqrt5)/2\}$. With the help of Lemma 2.1,
$$\align\sum_{k=0}^{p-1}\bi{2k}k\al^{2k}=&\al^{2p-2}\sum_{k=0}^{p-1}\bi{2k}k(-\al^{-1})^{2p-2-2k}
\\=&\al^{2p-2}\sum_{k=0}^{p-1}\bi{2k}k(\beta^2)^{p-1-k}
\\\eq&\al^{2p-2}\l(2u_p(\beta^2,\beta^2)-u_p(\beta^2-2,1)\r)\ (\mo\ p^2).
\endalign$$
By Lemma 2.3,
$$u_p(\beta^2-2,1)=u_p(\beta-1,1)=u_p(-\al,1)=U_{p-1}\l(\f{-\al}2\r)=U_{p-1}\l(\f{\al}2\r).$$
Note also that
$$\align u_p(\beta^2,\beta^2)=&\sum_{k=0}^{(p-1)/2}\bi{p-1-k}k(\beta^2)^{p-1-2k}(-\beta^2)^k
\\=&\beta^{p-1}\sum_{k=0}^{(p-1)/2}\bi{p-1-k}k(-1)^k\beta^{p-1-2k}=\beta^{p-1}U_{p-1}\l(\f{\beta}2\r).
\endalign$$
Therefore
$$\align&\al^{2p-2}\l(2u_p(\beta^2,\beta^2)-u_p(\beta^2-2,1)\r)
\\=&2\al^{p-1}(\al\beta)^{p-1}U_{p-1}\l(\f{\beta}2\r)-\al^{2p-2}U_{p-1}\l(\f{\al}2\r)
\\=&2\al^{p-1}U_{p-1}\l(\f{\beta}2\r)-\al^{2p-2}U_{p-1}\l(\f{\al}2\r).
\endalign$$

Observe that
$$U_{p-1}\l(\f{1+\sqrt5}4\r)=U_{p-1}\l(\cos\f{\pi}5\r)=\f{\sin(p\pi/5)}{\sin(\pi/5)}$$
and
$$U_{p-1}\l(\f{1-\sqrt5}4\r)=U_{p-1}\l(\f{\sqrt5-1}4\r)=U_{p-1}\l(\cos\f{2\pi}5\r)=\f{\sin(2p\pi/5)}{\sin(2\pi/5)}.$$

Combining the above we obtain the desired results. \qed

\proclaim{Lemma 2.4} For $n\eq\pm3\ (\mo\ 10)$, we have
$$\f{\sin(n\pi/5)}{\sin(\pi/5)}=\f{\pm 1\pm\sqrt5}2
\ \ \t{and}\ \ \f{\sin(2n\pi/5)}{\sin(2\pi/5)}=\f{\pm1\mp\sqrt5}2.$$
\endproclaim
\Proof. It suffice to note that
$$\f{\sin(3\pi/5)}{\sin(\pi/5)}=\f{\sin(2\pi/5)}{\sin(\pi/5)}=2\cos\f{\pi}5=\f{1+\sqrt5}2$$
and
$$\f{\sin(6\pi/5)}{\sin(2\pi/5)}=-\f{\sin(4\pi/5)}{\sin(2\pi/5)}=-2\cos\f{2\pi}5=\f{1-\sqrt5}2.$$
We are done. \qed

\proclaim{Lemma 2.5} Let
$p\not=2,5$ be a prime. If $p\eq\pm1\ (\mo\ 10)$, then
$$2F_{p-1}-F_{2p-2}\eq0\ (\mo\ p^2)\ \t{and}\ 2L_{p-1}-L_{2p-2}\eq2\ (\mo\ p^2).$$
If $p\eq\pm3\ (\mo\ 10)$, then
$$2F_{p-2}+F_{2p-1}\eq-2\ (\mo\ p^2)\ \t{and}\ 2L_{p-2}+L_{2p-1}\eq4\ (\mo\ p^2).$$
\endproclaim
\Proof. It is well known that for any $n\in\N$ we have
$$F_{2n}=F_nL_n\ \t{and}\ L_{2n}=L_n^2-2(-1)^n.$$
(See, e.g., [R].) By [SS] or [ST2], $L_{p-(\f p5)}\eq 2(\f p5)\ (\mo\ p^2)$.
Thus, if $p\eq\pm1\ (\mo\ 10)$ (i.e., $(\f p5)=1$), then
$$2F_{p-1}-F_{2p-2}=F_{p-1}(2-L_{p-1})\eq0\ (\mo\ p^2)$$
and
$$2L_{p-1}-L_{2p-2}=2L_{p-1}-(L_{p-1}^2-2)\eq2\times2-(2^2-2)=2\ (\mo\ p^2).$$

Now assume that $p\eq\pm3\ (\mo\ 10)$. Then
$$F_{p+1}=F_{p-(\f p5)}\eq0\ (\mo\ p),\ L_{p+1}=L_{p-(\f p5)}\eq2\l(\f p5\r)=-2\ (\mo\ p^2),$$
and
$$L_p=1+\f 52F_{p-(\f p5)}=1+\f 52 F_{p+1}\ (\mo\ p^2).$$
(Cf. [SS] and [ST2].) Thus
$$L_{p-2}=2L_p-L_{p+1}\eq 2+5F_{p+1}-(-2)=4+5F_{p+1}\ (\mo\ p^2)$$
and
$$\align L_{2p-1}=&L_{2p+2}-2L_{2p}=L_{p+1}^2-2-2(L_p^2+2)
\\\eq&(-2)^2-2-2\(\l(1+\f52F_{p+1}\r)^2+2\)\eq-4-10F_{p+1}\ (\mo\ p^2).
\endalign$$
It follows that
$$2F_{p-2}+L_{2p-1}\eq 8+10F_{p+1}-4-10F_{p+1}=4\ (\mo\ p^2).$$
Observe that
$$F_p=2F_{p+1}-L_p\eq 2F_{p+1}-\l(1+\f 52F_{p+1}\r)=-1-\f12F_{p+1}\ (\mo\ p^2)$$
and hence
$$F_{p-2}=2F_p-F_{p+1}\eq-2-2F_{p+1}\ (\mo\ p^2).$$
Note also that
$$\align F_{2p-1}=&F_{2p+2}-2F_{2p}=F_{p+1}L_{p+1}-2F_pL_p
\\\eq&-2F_{p+1}+2\l(1+\f12F_{p+1}\r)\l(1+\f 52F_{p+1}\r)
\\\eq&-2F_{p+1}+2(1+3F_{p+1})=4F_{p+1}+2\ (\mo\ p^2).
\endalign$$
Therefore
$$2F_{p-2}+F_{2p-1}\eq-4-4F_{p+1}+4F_{p+1}+2=-2\ (\mo\ p^2).$$

The proof of Lemma 2.5 is now complete. \qed

\proclaim{Lemma 2.6} Let $p$ be an odd prime.
Suppose that $\cos\theta$ is an algebraic $p$-adic integer. Then
$$\sum_{k=0}^{p-1}\bi{2k}k(2+2\cos\theta)^{p-1-k}=2(2\cos\theta+2)^{(p-1)/2}\ \f{\sin(p\theta/2)}{\sin(\theta/2)}
-\f{\sin(p\theta)}{\sin\theta}$$
and
$$\sum_{k=0}^{p-1}\bi{2k}k(2-2\cos\theta)^{p-1-k}=2(2\cos\theta-2)^{(p-1)/2}\ \f{\cos(p\theta/2)}{\cos(\theta/2)}
-\f{\sin(p\theta)}{\sin\theta}.$$
\endproclaim
\Proof. Note that if we replace $\theta$ in the first equality by $\pi-\theta$ we then get the second equality.
So it suffices to prove the first equality.

In view of Lemma 2.1, we have
$$\sum_{k=0}^{p-1}\bi{2k}k(2+2\cos\theta)^{p-1-k}
\eq2 u_p\l(4\cos^2\f{\theta}2,4\cos^2\f{\theta}2\r)-u_p(2\cos\theta,1)\ (\mo\ p^2).$$
With the help of Lemma 2.2,
$$\align &2 u_p\l(4\cos^2\f{\theta}2,4\cos^2\f{\theta}2\r)-u_p(2\cos\theta,1)
\\=&2\l(2\cos\f{\theta}2\r)^{p-1}U_{p-1}\l(\cos\f{\theta}2\r)-U_{p-1}(\cos\theta)
\\=&2(2+2\cos\theta)^{(p-1)/2}\ \f{\sin(p\theta/2)}{\sin(\theta/2)}-\f{\sin(p\theta)}{\sin\theta}.
\endalign$$
We are done. \qed

\proclaim{Lemma 2.7} Let $p$ be an odd prime. Then
$$U_{p-1}(x)\eq(-1)^{(p-1)/2}\sum_{k=0}^{(p-1)/2}\f{\bi{2k}k}{16^k}(4x^2)^k\ (\mo\ p^2).$$
\endproclaim
\Proof. Let $n=(p-1)/2$. By an observation of the author's brother Z. H. Sun, for $k=0,\ldots,n$ we have
$$\align\bi{n+k}{2k}=&\f{\prod_{0<j\ls k}(p^2-(2j-1)^2)}{4^k(2k)!}
\\\eq&\f{\prod_{0<j\ls k}(-(2j-1)^2)}{4^k(2k)!}
=\f{\bi{2k}k}{(-16)^k}\ (\mo\ p^2).
\endalign$$
Thus
$$\align U_{p-1}(x)=&\sum_{k=0}^n(-1)^k\bi{p-1-k}k(2x)^{p-1-2k}
\\=&\sum_{k=0}^n(-1)^k\bi{n+(n-k)}{n-(n-k)}(2x)^{2(n-k)}
\\=&\sum_{k=0}^n(-1)^{n-k}\bi{n+k}{n-k}(2x)^{2k}
\\\eq&(-1)^n\sum_{k=0}^n\f{\bi{2k}k}{16^k}(4x^2)^k\ (\mo\ p^2).
\endalign$$
This concludes the proof. \qed

\proclaim{Lemma 2.8} Let $p$ be an odd prime and set $p^*=(-1)^{(p-1)/2}p$.
Suppose that $\cos\theta$ is an algebraic $p$-adic integer. Then
$$\sum_{k=0}^{(p-1)/2}\f{\bi{2k}k}{16^k}\l((2+2\cos\theta)^k\pm(2-2\cos\theta)^k\r)
\eq2\f{\sin((p^*\pm1)\theta/2)}{\sin\theta}\ (\mo\ p^2).$$
\endproclaim
\Proof. Let $n=(p-1)/2$. Since
$$4\cos^2\f{\theta}2=2+2\cos\theta,$$
by Lemma 2.7 we have
$$(-1)^n\sum_{k=0}^{(p-1)/2}\f{\bi{2k}k}{16^k}(2+2\cos\theta)^k\eq U_{p-1}(\cos\theta)=\f{\sin(p\theta/2)}{\sin(\theta/2)}
\ (\mo\ p^2)$$
and also
$$\align &(-1)^n\sum_{k=0}^{(p-1)/2}\f{\bi{2k}k}{16^k}(2-2\cos\theta)^k
\\\eq& U_{p-1}(\cos(\pi-\theta))=\f{\sin(p(\pi-\theta)/2)}{\sin((\pi-\theta)/2)}
=(-1)^n\f{\cos(p\theta/2)}{\cos(\theta/2)}
\ (\mo\ p^2).
\endalign$$
Thus
$$\align &\sum_{k=0}^{(p-1)/2}\f{\bi{2k}k}{16^k}\l((2+2\cos\theta)^k\pm(2-2\cos\theta)^k\r)
\\\eq&(-1)^n\f{\sin(p\theta/2)}{\sin(\theta/2)}\pm\f{\cos(p\theta/2)}{\cos(\theta/2)}
=\f{\sin(p^*\theta/2)}{\sin(\theta/2)}\pm\f{\cos(p^*\theta/2)}{\cos(\theta/2)}
\\\eq&\f{\sin((p^*\pm1)\theta/2)}{(\sin\theta)/2}\ (\mo\ p^2).
\endalign$$
We are done. \qed

\heading{3. Proofs of Theorems 1.1-1.4}\endheading

\noindent{\it Proof of Theorem 1.1}. (i) Let us first prove the first and the second congruences
in Theorem 1.1. Since $L_{2k}+F_{2k}=2F_{2k+1}$ for any $k\in\N$,
we only need to show the first congruence
and the following one:
$$\sum_{k=0}^{p-1}L_{2k}\bi{2k}k\eq(-1)^{\lfloor p/5\rfloor}\l(3\l(\f p5\r)-1\r)\ (\mo\ p^2).$$

 Set
$$\al=\f{1+\sqrt5}2\ \ \t{and}\ \ \beta=\f{1-\sqrt5}2.$$
By Lemma 2.3 we have
$$\align&\sum_{k=0}^{p-1}\bi{2k}kL_{2k}=\sum_{k=0}^{p-1}\bi{2k}k(\al^{2k}+\beta^{2k})
\\=&2\al^{p-1}\f{\sin(2p\pi/5)}{\sin(2\pi/5)}-\al^{2p-2}\f{\sin(p\pi/5)}{\sin(\pi/5)}
\\&+2\beta^{p-1}\f{\sin(p\pi/5)}{\sin(\pi/5)}-\beta^{2p-2}\f{\sin(2p\pi/5)}{\sin(2\pi/5)}.
\endalign$$
Thus, if $p\eq\pm1\ (\mo\ 10)$ then
$$\sum_{k=0}^{p-1}\bi{2k}kL_{2k}=\pm(2\al^{p-1}-\al^{2p-2}+2\beta^{p-1}-\beta^{2p-2})=\pm(2L_{p-1}-L_{2p-2}).$$
When $p\eq\pm3\ (\mo\ 10)$, by Lemma 2.4 and the above we have
$$\align\sum_{k=0}^{p-1}\bi{2k}kL_{2k}
=&2\al^{p-1}(\pm\beta)-\al^{2p-2}(\pm\al)+2\beta^{p-1}(\pm\al)-\beta^{2p-2}(\pm\beta)
\\=&\pm2\al\beta(\al^{p-2}+\beta^{p-2})\mp(\al^{2p-1}+\beta^{2p-1})=\mp2L_{p-2}\mp L_{2p-1}.
\endalign$$

In light of Lemma 2.3, we also have
$$\align&\sqrt5\sum_{k=0}^{p-1}\bi{2k}kF_{2k}=\sum_{k=0}^{p-1}\bi{2k}k(\al^{2k}-\beta^{2k})
\\=&2\al^{p-1}\f{\sin(2p\pi/5)}{\sin(2\pi/5)}-\al^{2p-2}\f{\sin(p\pi/5)}{\sin(\pi/5)}
\\&-2\beta^{p-1}\f{\sin(p\pi/5)}{\sin(\pi/5)}+\beta^{2p-2}\f{\sin(2p\pi/5)}{\sin(2\pi/5)}.
\endalign$$
Thus, if $p\eq\pm1\ (\mo\ 10)$ then
$$\align\sqrt5\sum_{k=0}^{p-1}\bi{2k}kF_{2k}=&\pm(2\al^{p-1}-\al^{2p-2}-2\beta^{p-1}+\beta^{2p-2})
\\=&\pm\sqrt5(2F_{p-1}-F_{2p-2}).
\endalign$$
When $p\eq\pm3\ (\mo\ 10)$, by Lemma 2.4 and the above we have
$$\align\sqrt5\sum_{k=0}^{p-1}\bi{2k}kF_{2k}
=&2\al^{p-1}(\pm\beta)-\al^{2p-2}(\pm\al)-2\beta^{p-1}(\pm\al)+\beta^{2p-2}(\pm\beta)
\\=&\pm2\al\beta(\al^{p-2}-\beta^{p-2})\mp(\al^{2p-1}-\beta^{2p-1})
\\=&\mp\sqrt5(2F_{p-2}+F_{2p-1}).
\endalign$$

In view of the above and Lemma 2.5, we have proved the first and the second congruences in Theorem 1.1.

\medskip

(ii) Recall that
$$\cos\f{2\pi}5=\f{\sqrt5-1}4\ \t{and}\ \cos\f{\pi}5=\f{\sqrt5+1}4.$$
Applying Lemma 2.7 we get
$$(-1)^n\sum_{k=0}^n\f{\bi{2k}k}{16^k}\l(\f{1-\sqrt5}2\r)^{2k}\eq U_{p-1}\l(\cos\f{2\pi}5\r)=\f{\sin(2p\pi/5)}{\sin(2\pi/5)}
\ (\mo\ p^2)$$
and
$$(-1)^n\sum_{k=0}^n\f{\bi{2k}k}{16^k}\l(\f{1+\sqrt5}2\r)^{2k}\eq U_{p-1}\l(\cos\f{2\pi}5\r)
=\f{\sin(p\pi/5)}{\sin(\pi/5)}\ (\mo\ p^2).$$
Combining this with Lemma 2.4 we can easily deduce the third and the fourth congruences in Theorem 1.1.
(Note that $2F_{2k+1}=F_{2k}+L_{2k}$ for $k\in\N$.)
\medskip

So far we have completed the proof of Theorem 1.1.

\medskip
\noindent{\it Proof of Theorem 1.2}. The proof is similar to that of Theorem 1.1. Observe that
$$4\cos^2\f{\pi}{10}=2+2\cos\f{\pi}5=\f{5+\sqrt5}2
\ \t{and}\ 4\cos^2\f{3\pi}{10}=2+2\cos\f{3\pi}5=\f{5-\sqrt5}2.$$
So we may employ the results on $F_{(p\pm1)/2}$ and $L_{(p\pm1)/2}$ modulo $p^2$ in [SS]
to get the four congruences with helps of Lemmas 2.6 and 2.8.
\qed

\medskip
\noindent{\it Proof of Theorem 1.3}. To get the desired congruences we may apply Lemmas 2.6 and 2.8 with $\theta=\pi/6$
and use the results on $S_{(p\pm1)/2}$ and $T_{(p\pm1)/2}$ modulo $p^2$ in [S02]. \qed

\medskip
\noindent{\it Proof of Theorem 1.4}. Apply Lemmas 2.6 and 2.8 with $\theta=\pi/4$ and use a result in [Su]. \qed

\heading{4. Some open conjectures}\endheading

In this section we formulate several open conjectures.

\proclaim{Conjecture 4.1} For any $n\in\Z^+$ we have
$$\f{(-1)^{\lfloor n/5\rfloor-1}}{(2n+1)n^2\bi{2n}n}\sum_{k=0}^{n-1}F_{2k+1}\bi{2k}k\eq
\cases6\ (\mo\ 25)&\t{if}\ n\eq0\ (\mo\ 5),\\4\ (\mo\ 25)&\t{if}\ n\eq1\ (\mo\ 5),
\\1\ (\mo\ 25)&\t{if}\ n\eq2,4\ (\mo\ 5),\\9\ (\mo\ 25)&\t{if}\ n\eq3\ (\mo\ 5).\endcases$$
Also, if $a,b\in\Z^+$ and $a\gs b$ then the sum
$$\f1{5^{2a}}\sum_{k=0}^{5^a-1}F_{2k+1}\bi{2k}k$$
modulo $5^b$ only depends on $b$.
\endproclaim

Recall that the usual $q$-analogue of $n\in\N$ is given by
$$[n]_q=\f{1-q^n}{1-q}=\sum_{0\ls k<n}q^k$$
which tends to $n$ as $q\to1$.
For any $n,k\in\N$ with $n\gs k$,
$$\M nk=\f{\prod_{0<r\ls n}[r]_q}{(\prod_{0<s\ls k}[s]_q)
(\prod_{0<t\ls n-k}[t]_q)}$$
is a natural extension of the usual binomial coefficient
$\bi nk$. A $q$-analogue of Fibonacci numbers introduced by I. Schur [Sc]
is defined as follows:
$$F_0(q)=0,\ F_1(q)=1, \t{and}\ F_{n+1}(q)=F_n(q)+q^nF_{n-1}(q)\ (n=1,2,3,\ldots).$$

\proclaim{Conjecture 4.2} Let $a$ and $m$ be positive integers.
 Then we have the following congruence in the ring $\Z[q]$:
$$\sum_{k=0}^{5^am-1}q^{-2k(k+1)}\M{2k}kF_{2k+1}(q)\eq0\ (\mo\ [5^a]_q^2).$$
\endproclaim

\proclaim{Conjecture 4.3} For any $n\in\Z^+$ we have
$$\f{(-1)^{n-1}}{n^2(n+1)\bi{2n}n}\sum_{k=0}^{n-1}S_{k+1}\bi{2k}k\eq
\cases1\ (\mo\ 9)&\t{if}\ n\eq0,2\ (\mo\ 9),\\4\ (\mo\ 9)&\t{if}\ n\eq5,6\ (\mo\ 9),
\\-2\ (\mo\ 9)&\t{otherwise}.\endcases$$
Also, if $a,b\in\Z^+$ and $a\gs b-1$ then the sum
$$\f1{3^{2a}}\sum_{k=0}^{3^a-1}S_{k+1}\bi{2k}k$$
modulo $3^b$ only depends on $b$.
\endproclaim

\proclaim{Conjecture 4.4} Let $p>3$ be a prime. Then
$$\sum_{k=0}^{p-1}\bi{p-1}k\bi{2k}k((-1)^k-(-3)^{-k})\eq\l(\f p3\r)(3^{p-1}-1)\ (\mo\ p^3).$$
\endproclaim
\Remark. The congruence mod $p^2$ follows from [S09b].

\proclaim{Conjecture 4.5} Let $p$ be a prime with $p\eq\pm1\ (\mo\ 12)$. Then
$$\sum_{k=0}^{p-1}\bi{p-1}k\bi{2k}k(-1)^kS_k\eq(-1)^{(p-1)/2}S_{p-1}\ (\mo\ p^3).$$
\endproclaim
\Remark. The author can prove the congruence mod $p^2$.

\medskip

\proclaim{Conjecture 4.6} Let $p$ be a prime with $p\eq\pm1\ (\mo\ 8)$. Then
$$\sum_{k=0}^{p-1}\bi{p-1}k\bi{2k}k\f{u_k(4,2)}{(-2)^k}\eq(-1)^{(p-1)/2}u_{p-1}(4,2)\ (\mo\ p^3).$$
\endproclaim
\Remark. The author has proved the  congruence mod $p^2$.

\enddocument